\documentclass[11pt,bezier,epsf]{article}
\usepackage{amsmath,amssymb,amsfonts}
   \usepackage{amsmath}
   \usepackage{amsthm}
   \usepackage{amssymb}
   \usepackage{latexsym}

\newtheorem{thm}{Theorem}[section]

\newtheorem{lem}[thm]{Lemma}

\newtheorem{defn}[thm]{Definition}

\newtheorem{exam}[thm]{Example}
\numberwithin{equation}{section}

\newcommand{\Rset}{\mathbb R}

\newcommand{\RF}{{\mathbb R}_{\mathcal F}}
\newcommand{\UP}{\underline{u}}
\newcommand{\UB}{\overline{u}}

\newcommand{\To}{\rightarrow}

\begin{document}
\title{A Note on "Global solutions for nonlinear fuzzy  fractional integral and integrodifferential equations"}
\author{%
   R. Alikhani$^{}$\footnote{Corresponding Author. Tel.: +1(321)674-7213;
Fax.:+1(321)674-7412.\
  \ \
\newline
   Email addresses: alikhani@tabrizu.ac.ir (R. Alikhani), fbahram@tabrizu.ac.ir (F. Bahrami).}, F. Bahrami\\
{\small \em Faculty of Mathematical Sciences, University of Tabriz,
Tabriz, Iran.}} \maketitle
\begin{abstract}
The authors in \cite{alikhani} have given two examples to illustrate
their results in which they have been eliminated the technical
details. However, the authors in \cite{salahshur} claimed that the
examples are incorrect. In fact they conjectured that the authors in
\cite{alikhani} have employed the incorrect statement $x-x=0$ for
$x\in \RF\setminus \Rset$ to construct the examples. Here we intend
to observe that the basic method used in \cite{alikhani} to prove
the validity of the examples is the well-known L-U representation of
a fuzzy-number valued function. In this sense, we will make use of
the $\alpha$-level sets and show that the examples are correct.
\end{abstract}
\section{Introduction}
The authors in \cite{alikhani}   have considered  the following
fuzzy fractional integrodifferential equations (FFIDEs) involving
Riemann-Liouville
 derivatives of  order $0 < q < 1$,
\begin{eqnarray}\label{fractional differential equation}
\label{maineq} &&D^qu(t)=f(t,u, Tu), \ \ \ \ \ \forall\ t\in
J,\nonumber\\
&&\lim_{t\to 0^+} t^{1-q}u(t)=u_0,
\end{eqnarray}
where $J=(0,b]$, $u_0\in \RF$,  $f\in C(J\times \Rset_\mathcal F^2,
\Rset_\mathcal F)$,
 $D^q u$ denotes fuzzy derivative of fractional order $q$ introduced in
 \cite{alikhani}, and
\begin{equation}\label{integral part}
 (Tu)(t)= \int_0^t k(t,s)u(s)ds,
\end{equation}
 where $k\in C(I,\Rset_+)$, $I=\{ (t,s)\in \bar{J}\times \bar{J}:t\geq s\} $, $\Rset_+=[0,+\infty)$ and $0<q<1$.
 To prove the existence theorems, they applyed concepts  of the upper
 and lower solutions for more generic-following fuzzy fractional integral equation
\begin{equation*}\label{fractional integral equation}
u(t)=g(t)+\frac{1}{\Gamma (q)}\int_0 ^t
(t-s)^{q-1}f(s,u(s),(Tu)(s))ds,\ \ \ \ t\in J,
\end{equation*}
where $J=(0,b]$,  $g\in C(J, \Rset_\mathcal F)$,  $f\in C(J\times
\Rset_\mathcal F^2, \Rset_\mathcal F)$.
 They combined these concepts with the monotone iterative technique and proved the existence  of solutions for Problem (\ref{fractional differential
 equation}).\\
  They have given two examples
to illustrate  the method in which they gave the upper and lower
solutions and the exact solution corresponded with examples.
However, while the proof of details was almost clear, then they
eliminated them. On the other hand the authors in \cite{salahshur}
have asserted that there exist some errors in the given examples in
\cite{alikhani}, and the given exact solutions are incorrect. In
fact they supposed the authors in \cite{alikhani} have employed the
incorrect expression $x-x=\hat{0}$  for $x\in \RF\setminus \Rset$ to
prove the details of their assertion. In this study  we recall the
definition of the solution for Problem (\ref{fractional differential
equation}) in \cite{alikhani}, and we prove the exact solutions by
verify the corresponding problems using $\alpha$-level sets not by
applying the incorrect expression $x-x=\hat{0}$. Also the authors in
\cite{salahshur} have asserted these solutions do not lie between
the lower and upper solutions and have given some examples to prove
it. In this paper we will state that according to Theorem 4.5-4.6 in
\cite{alikhani}, it is not necessary that all of the solutions lie
between the lower and upper solutions, it is true if this happens
just for some solutions not for all.
\section{Preliminaries}
In this section we recall a few known results that are needed in our
discussion. \\

$\RF$ denotes the space of fuzzy numbers on $\Rset$.
 For $0<\alpha\leq 1$, $\alpha$-level set of $x\in
\Rset_{\mathcal F}$ is defined by $[x]^{\alpha}=\{t \in \Rset \mid
x(t)\geq\alpha\}$ and $[x]^{0}=\overline{\{t \in \Rset \mid x(t)>
0\}}$. For any $\alpha \in [0,1]$, $[x]^{\alpha}$ is a bounded
closed interval, we denote
  $[x]^{\alpha}=[x_{l\alpha },x_{r\alpha }]$. Also, we define
  $\hat{0}\in\RF$ as $\hat{0}(t)=1$ if $t=0$ and $\hat{0}(t)=0$ if $t\neq
  0$. We denote a triangular fuzzy
 number as $x=(a,b,c)$, where $a$ and $c$ are endpoints of the $0$-level set and $1$-level set
 $=\{b\}$.\\
 For the proof of the validity of the examples in \cite{alikhani},
 we need to utilize $\alpha$-level sets concepts. In this sense, we
 will make use of the following well-known lemma.
 \begin{lem}(See e.g. \cite{Goetschel and Voxman}).\label{condition of fuzzy
number} Assume  $x_l:[0,1] \To \Rset $ and $x_r: [0,1]\To \Rset$
satisfy the conditions:
\item (i) $x_l $ is a bounded increasing function.
\item (ii) $x_r$ is a bounded decreasing function.
\item (iii) $x_l(1)\leq x_r(1)$.
\item (iv) For $0<k\leq 1$, $\lim_{\alpha \To k^-} x_l(\alpha)=x_l(k)$
and $\lim_{\alpha \To k^-} x_r(\alpha)=x_r(k)$.
\item (v) $\lim_{\alpha \To 0^+} x_l(\alpha)=x_l(0)$ and $\lim_{\alpha \To 0^+}
x_r(\alpha)=x_r(0)$.\\
 Then $[x_l(\alpha), x_r(\alpha)]$ is the
parametric form of a fuzzy number.
\end{lem}
  \begin{lem}(See e.g. \cite{Anastassiou}.)\label{properties of fuzzy number}\\
(i) If we denote $\hat{0}=\chi_{\{0\}}$, then
 $\hat{0}\in\Rset_{\mathcal F}$ is neutral element with respect to
+, i.e. $u+\hat{0}=\hat{0}+u=u,$ for all $u\in
\Rset_{\mathcal F}.$\\
(ii) With respect to $\hat{0},$\ none of\ $u\in \Rset_{\mathcal
F}\setminus \Rset,$\text{ has opposite in} $\Rset_{\mathcal F}$(\ with respect to +).\\
(iii) For any $a,b\in \Rset$ with $a,b\geq 0$ or $a,b\leq 0$ and any
$u\in \Rset_{\mathcal F},$ we have $(a+b)\cdot u=a\cdot u +b\cdot
u.$ For general $a,b\in \Rset,$ the above property does
not hold.\\
(iv) For any $\lambda \in \Rset$ and any $u,v \in \Rset_{\mathcal
F},$ we have $\lambda\cdot(u+v)=\lambda \cdot u+\lambda\cdot v.$\\
(v) For any $\lambda ,\mu\in \Rset$ and any $u\in \Rset_{\mathcal
F},$ we have $\lambda\cdot(\mu\cdot u)=(\lambda \mu)\cdot u.$
\end{lem}
 The space of continuous fuzzy
functions is denoted by $C(\bar{J},\Rset_\mathcal{
 F})$.
Moreover, let $0<r<1$. We consider the following space
\[C_{r}(\bar{J},\RF)=\{g\in C(J,\RF): t^{r}g\in C(\bar{J},\RF)\}.\]
The  definition  of solution for the problem (\ref{fractional
differential equation}) has been recalled from \cite{alikhani} as
follows.
\begin{defn}\label{definition of solution}
We say that $u\in C_{1-q}(\bar{J},\RF)$ is a solution for the
problem (\ref{fractional differential
  equation}), if it satisfies Problem (\ref{fractional differential
  equation}).
\end{defn}
In the following examples we intend to state the details of Examples
1 and 2 in \cite{alikhani}  and also to show  what is wrong in
Example 1.1-1.2 \cite{salahshur}, respectively.
\begin{exam}
Consider  the  linear fuzzy initial value problem
\begin{eqnarray}\label{exam1}
&&D^qu(t)=(\frac{t^{-q}}{\Gamma(1-q)}-t)u(t)+\int_0^t u(s)ds, \ \ \ \ t\in(0,b],\nonumber\\
&&\lim_{t\to 0^+}t^{1-q}u(t)=\hat{0},
\end{eqnarray}
where $0<q<1$ and $b=\sqrt[1+q]{\frac{1}{\Gamma(1-q)}}$.
$\UP_1=\hat{0}$ and $\UB_1=t^q$ are the lower and upper solutions of
Problem (\ref{exam1}), respectively, and $u(t)=c $ for $c\in \RF$ is
an exact solution as stated in \cite{alikhani}. Here in details we
prove  $u(t)=c$, for all $c\in \RF$ is a solution for Problem
(\ref{exam1}) based on Definition \ref{definition of solution}. To
this end we have to prove
\[D^q c=(\frac{c^{-q}}{\Gamma(1-q)}-t)c+\int_0^t cds, \ \ \ \ t\in(0,b],\]
i.e.
\begin{equation}\label{expression}
\frac{t^{-q}}{\Gamma(1-q)}c=(\frac{t^{-q}}{\Gamma(1-q)}-t)c+t c, \ \
\ \ t\in(0,b].
\end{equation}
We consider $\alpha$-level set of $c$ as
$[c]^{\alpha}=[c_{l\alpha},c_{r\alpha}]$. Since $0<t\leq b$, the
sign of $\frac{t^{-q}}{\Gamma(1-q)}-t$ is positive. Thus we can
write the right-hand (\ref{expression}) as follows:
\begin{eqnarray}
(\frac{t^{-q}}{\Gamma(1-q)}-t)[c_{l\alpha},c_{r\alpha}]+t
[c_{l\alpha},c_{r\alpha}]&=&[(\frac{t^{-q}}{\Gamma(1-q)}-t)c_{l\alpha}+tc_{l\alpha},(\frac{t^{-q}}{\Gamma(1-q)}-t)c_{r\alpha}+tc_{r\alpha}]\nonumber\\
&=&[\frac{t^{-q}}{\Gamma(1-q)}c_{l\alpha},\frac{t^{-q}}{\Gamma(1-q)}c_{r\alpha}]\nonumber\\
&=&\frac{t^{-q}}{\Gamma(1-q)}[c_{l\alpha},c_{r\alpha}]\nonumber
\end{eqnarray}
This proves the equality of (\ref{expression}). It then follows that
$c$ is an exact solution of Problem (\ref{exam1}). The authors in
\cite{salahshur} have guessed that in the proof of equality
(\ref{expression}) in Example 1 in \cite{alikhani}, $tc-tc=0$ has
been used. It sounds like the authors in \cite{salahshur}  have
employed the following idea to prove Eq. (\ref{expression})
\[(\frac{t^{-q}}{\Gamma(1-q)}-t)c+t c= \frac{t^{-q}}{\Gamma(1-q)}c-tc+tc=\frac{t^{-q}}{\Gamma(1-q)}c.\]
By virtue of  Lemma \ref{properties of fuzzy number}(ii)-(iii), it
is not possible. Since
 $\UP_1=\hat{0}$ and $\UB_1=t^q$, all  the conditions in  Theorem 4.5 in \cite{alikhani} are fulfilled. Then there
exists a solution of the problem (\ref{exam1}), $u$,  so that
$\hat{0}\leq u(t)\leq t^q$ for  $t\in (0,b]$,  $q\in (0,1)$. Here
this happens for $u(t)=c=\hat{0}$. The authors in \cite{salahshur}
have misunderstood that  this should be happen for all $c\in \RF$,
and they have given an example as a contraction example.
\end{exam}
\begin{exam}
Consider  the  linear fuzzy initial value problem
\begin{eqnarray}\label{exam2}
&&D^qu(t)=\frac{c}{\Gamma(1-q)}({t^{-q}}-1-t^{q-1})+\frac{1}{\Gamma(1-q)}u(t), \ \ \ \ t\in(0,0.32],\nonumber\\
&&\lim_{t\to 0^+}t^{1-q}u(t)=c,
\end{eqnarray}
where $0.58<q\leq0.88$ and $c\in \RF$.  $\UP_1=ct^{q-1}$ and
$\UB_1=10ct^{q-1}$ are the lower and upper solutions of  Problem
(\ref{exam2}), respectively, and $u(t)=c+ct^{q-1} $  is an exact
solution as stated in \cite{alikhani}. Here  we prove in details
that $u(t)=c+ct^{q-1} $, for all $c\in \RF$ is a solution for
Problem (\ref{exam2}) based on  Definition \ref{definition of
solution}. To this end, we have to prove that
\[D^q (c+ct^{q-1})=(\frac{c}{\Gamma(1-q)}(t^{-q}-1-t^{q-1})+\frac{c+ct^{q-1}}{\Gamma(1-q)}, \ \ \ \ t\in(0,0.32],\]
i.e.
\begin{equation}\label{expression2}
\frac{ct^{-q}}{\Gamma(1-q)}=\frac{c}{\Gamma(1-q)}(t^{-q}-1-t^{q-1})+\frac{c+ct^{q-1}}{\Gamma(1-q)},
\ \ \ \ t\in(0,0.32].
\end{equation}
We consider $\alpha$-level set of $c$ as
$[c]^{\alpha}=[c_{l\alpha},c_{r\alpha}]$. Since $ t\in(0,0.32]$, the
sign of $t^{-q}-1-t^{q-1}$ is positive. Thus, we can prove the
equality below in a similar way to Example(\ref{exam1}).
\begin{equation*}
\frac{[c_{l\alpha},c_{r\alpha}]t^{-q}}{\Gamma(1-q)}=\frac{[c_{l\alpha},c_{r\alpha}]}{\Gamma(1-q)}(t^{-q}-1-t^{q-1})+\frac{[c_{l\alpha},c_{r\alpha}]+[c_{l\alpha},c_{r\alpha}]t^{q-1}}{\Gamma(1-q)},
\ \ \ \ t\in(0,0.32].
\end{equation*}

This proves Eq.(\ref{expression2}), and then $u(t)=c+ct^{q-1} $ is
the exact solution of Problem (\ref{exam2}). Since all the
conditions of Theorem 4.5 in \cite{alikhani} are fulfilled, there
exists a solution $u$ for this problem  so that $ct^{q-1}\leq
u(t)\leq 10ct^{q-1}$ for all $t\in (0,0.32]$ and $0.58<q\leq0.88$.
Here this happens  for $u(t)=c+ct^{q-1}$ where $c\geq\hat{0}$ not
for all $c\in \RF$.
\end{exam}

\end{document}